\documentclass[11pt]{article} 
 \usepackage{amssymb}
\usepackage{amsthm}
 \usepackage{graphicx}
\usepackage{indentfirst}

\begin{document}
\begin{center}
\LARGE\noindent\textbf{On Hamiltonian bypasses in digraphs and bipartite digraphs}\\

\end{center} 

\begin{center}
\noindent\textbf{Samvel Kh. Darbinyan}\\

Institute for Informatics and Automation Problems,

Armenian National Academy of Sciences

E-mail: samdarbin@iiap.sci.am \\

\end{center}

\textbf{Abstract}

A Hamiltonian path in a digraph $D$ in which the initial vertex dominates the terminal vertex is called a Hamiltonian bypass.
Let $D$ be a 2-strong digraph of order $p\geq 3$ and let $z$ be some vertex of $D$. Suppose that every vertex of $D$ other than $z$ has degree at least $p$.  We introduce and study a conjecture which claims that there exists a smallest integer $k$ such that if $d(z)\geq k$, then $D$ contains a Hamiltonian bypass.
In this paper, we prove: 
(i)  If $D$ is Hamiltonian or $z$ has a degree greater than $(p-1)/3$,  
 then $D$ contains a Hamiltonian bypass.
(ii)  If a strong balanced bipartite digraph $B$ of order $2a\geq 6$ satisfies the condition that $d^+(u)+d^-(v)\geq a+1$ for all vertices $u$ and $v$ from different partite sets such that $B$ does not contain the arc $uv$, then $B$ contains a Hamiltonian bypass. Furthermore, the lower bound $a+1$  is sharp. The first result improves a result of Benhocine (J. of Graph Theory, 8, 1984) and a result of the author (Math. Problems of Computer Science,  
 54, 2020)

We also suggest some conjectures and problems.

 
  \textbf{AMS} subject classifications: 05C20, 05C38, 05C45\\

\section {Introduction} 

 We shall assume that the reader is familiar with the standard terminology on digraphs and  for the terminology and notation not defined in this paper,   the reader is referred to \cite{[2]}. 
In this paper, we consider finite digraphs without loops and multiple arcs, which may contain opposite arcs with the same endvertices, i.e., a cycle of length two. Each cycle and path is assumed to be simple and directed. 
A cycle (path) in a digraph $D$ which passes through all the vertices of $D$ is called {\it Hamiltonian}. 
A digraph containing a Hamiltonian cycle is called a {\it Hamiltonian digraph}. One of the fundamental and most studied problems in digraph theory is to find sufficient conditions for a digraph to contain a Hamiltonian  path of a certain type. In a digraph $D$, a Hamiltonian path  is called a {\it Hamiltonian bypass} if its  initial vertex dominates  its terminal vertex. 
 There are many sufficient conditions for the existence of a Hamiltonian cycle in  digraphs (see, e.g., \cite{[3]}, \cite{[12]},  \cite{[13]}, \cite{[15]}, \cite{[18]}, \cite{[19]}, \cite{[20]}, \cite{[22]}). It is natural to consider  an analogous problem for  the existence of  a Hamiltonian bypass.

 It was proved in  \cite{[4]}, \cite{[5]}, \cite{[8]}, \cite{[9]}, \cite{[10]}, \cite{[14]}, that a number of   sufficient conditions for a digraph to be Hamiltonian are also sufficient  for a digraph to contain a Hamiltonian bypass (except for  some exceptional digraphs, which are characterised). These include  Theorems 1.1-1.4.     To formulate these theorems, we need the following definitions.

\textbf{Definition 1.}  Let $D_0$ be an arbitrary digraph of order $p\geq 3$ with $p$  odd, such that $V(D_0)=F\cup B$, where $F\cap B=\emptyset$, $F$ is an independent set with $(p+1)/2$ vertices, $B$ is a set of $(p-1)/2$ vertices that  induce an arbitrary subdigraph. Furthermore,  $D_0$ contains all the possible arcs between $F$ and $B$.
 
By $K^*_m$ we denote the complete digraph of order $m$.

\textbf{Definition 2.}  For any $k\in [1, p-2]$ let $D_{p-k, k}$ denote a digraph of order $p\geq 3$, obtained from $K^*_{p-k}$ and $K^*_{k+1}$ by identifying a vertex of the former with a vertex of the latter.
 
\textbf{Definition 3.}  By $T_5$ we denote a tournament of order 5 with vertex set $\{x_1, x_2,$ $ x_3, x_4, y\}$ and arc set $\{x_ix_{i+1}\,| \, i\in [1,3]\}\cup \{x_4x_1, x_1y, x_3y, yx_2, yx_4,x_1x_3, x_2x_4\}$.

 \textbf{Theorem 1.1} (Benhocine \cite{[4]}). {\it Let $D$ be a 2-strong digraph of order $p$ with minimum degree at least $p-1$. Then $D$ contains a Hamiltonian bypass, unless $D$ is isomorphic to a digraph of type $D_0$.}

 \textbf{Theorem 1.2} (Benhocine \cite{[4]}). {\it Every digraph $D$ of order $p\geq 3$ and minimum degree at least $p$ contains a Hamiltonian bypass}.
 
\textbf{Theorem 1.3} (Darbinyan \cite{[8]}). {\it Let $D$ be a strong digraph of order $p\geq 3$. Suppose that $d(x)+d(y)\geq 2p-2$ 
 for every pair of non-adjacent vertices $x$, $y$ of $V(D)$.  Then $D$ contains a Hamiltonian bypass, unless $D$ is isomorphic to a digraph of the set  $D_0\cup \{D_{p-k,k},T_5,C_3\}$,  where $C_3$ is the directed cycle of length three.}

\textbf{Theorem 1.4} (Darbinyan \cite{[10]}). {\it Let $D$ be a strong digraph of order $p\geq 4$. Suppose that  
 for each triple of vertices $x$, $y$, $z$ such that  $x$ and $y$ are non-adjacent vertices: If there is no arc from  $x$ to  $z$, then  $d(x)+d(y) +d^+(x)+d^-(z)\geq 3p-2$.  If there is no arc from  $z$ to  $x$, then  $d(x)+d(y) +d^-(x)+d^+(z)\geq 3p-2$. Then $D$ contains a Hamiltonian bypass, unless $D$ is isomorphic to the tournament $T_5$.}\\

We denote by  $H(n)$ the digraph of order $n\geq 8$ with vertex set 
$$
V(H(n)):=\{x_0,x_1,x_2,\ldots , x_{n-4}, y_1, y_2, y_3\} \quad \hbox{and arc set}
$$
$$
A(H(n)):=\{y_iy_j \, | \, i\not= j\} \cup \{x_ix_{i+1}\, | \, 0\leq i\leq n-4\} $$ $$\cup \{y_ix_j \, | \, 1\leq i\leq 3, 1\leq j\leq n-6\}
  \cup \, \{x_ix_j \, | \, 1\leq j< i\leq n-4\}$$ $$ \cup \{x_{n-4}y_i, x_{n-6}y_i \, | \, 1\leq i\leq 3\} \cup \{x_ix_{n-5}\, | \, 1\leq i\leq n-7\}$$
  $$
  \cup \, \{x_0x_{n-5},x_{n-5}x_0, x_{n-4}x_0, x_{n-6}x_{n-4}\}.
$$
It is easy to check that  $n-1$ vertices  of $H(n)$ have degrees  at least $n$ and its minimum degree equal to four. 
In \cite{[11]}, the author proved  that $H(n)$ is 2-connected and is not Hamiltonian. 

Let $D$ be a 2-strong digraph of order $p\geq 3$ in which $p-1$ vertices have degree at least $p$. In \cite{[9]}, the author proved that if $D$ is Hamiltonian or its the minimum degree  is  more than $0.4(p-1)$, then $D$ has a Hamiltonian bypass. In this paper, we improve this result, by proving the following theorem.

\textbf{Theorem 1.5.} {\it Let $D$ be a 2-strong digraph of order $p\geq 3$ and let  $z$ be some vertex of $D$. Suppose that $d(x)\geq p$ for each  vertex $x\in  V(D)\setminus \{z\}$.  If $D$ is Hamiltonian or $d(z)>(p-1)/3$,  then $D$ has a Hamiltonian bypass.}\\

Let $B(X,Y)$ be a balanced bipartite digraph with partite sets $X=\{x_1,x_2,\ldots , x_a\}$ and $Y=\{y_1,y_2,\ldots , y_a\}$. For an integer $l\geq 0$, we will say that $B(X,Y)$ satisfies the condition $A_l$ if $d^+(u)+d^-(v)\geq a+l$ for all $u$ and $v$ from different partite sets such that $uv\notin A(B(X,Y))$.

In \cite{[1]}, Adamus and Adamus proved: If $a\geq 2$ and a balanced bipartite digraph satisfies the condition $A_2$,  then it is Hamiltonian. In \cite{[1]}, the authors presented  an example of balanced bipartite digraph $B_6$ of order 6 with partite sets  $X=\{x_1,x_2,x_3\}$ and $Y=\{y_1,y_2,y_3\}$
(see Figure 1(b), in Figure 1(b) an undirected edge represents two directed arcs of opposite directions). It is easy to check that the digraph $B_6$ satisfies the condition $A_1$, but is not Hamiltonian. 
Wang \cite{[21]} proved that if a  balanced bipartite digraph  of order $2a\geq 6$ other than $B_6$  satisfies the condition $A_1$,  
is also  Hamiltonian. Note that $D(p,2)=[x_1y_1; x_1y_3x_2y_2x_3y_1]$ is a Hamiltonian bypass in $B_6$.

T\label{Fig.1}

 \begin{figure*}[!h]
\centerline{\includegraphics[width=10cm]{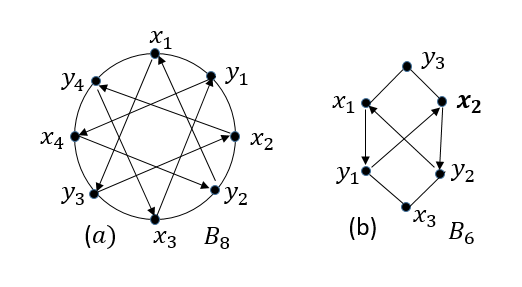}}
\caption{(a) The balanced bipartite digraph $B_8$ has no Hamiltonian bypass; (b) The balanced bipartite digraph $B_6$ is not Hamiltonian.}
\label{Fig.1}
 \end{figure*}

 In this paper, we also prove the following theorem for strong balanced bipartite digraphs.
 
 \textbf{Theorem 1.6.}
 {\it If a strong balanced bipartite digraph $B$ of order $2a\geq 6$ satisfies the condition that $d^+(u)+d^-(v)\geq a+1$ for all $u$ and $v$ from different partite sets such that $uv\notin A(B(X,Y))$, then it contains a Hamiltonian bypass}.\\

We shall prove a slightly more general result  (see, Theorem 4.1).
 
\section {Terminology and notation}

  In this paper, we consider finite digraphs (directed graphs) without loops and multiple arcs.  For a digraph $D$, we denote
  by $V(D)$ the vertex set of $D$ and by  $A(D)$ the set of arcs in $D$. The {\it order} of $D$ is the number
  of its vertices. Let $x$, $y$ be distinct vertices in $D$.
 The arc of a digraph $D$ directed from
   $x$ to $y$ is denoted by $xy$ (we say that $x$ dominates $y$). When $D$ is known from the context, sometimes when $xy\in A(D)$ or $xy\notin A(D)$, we will write $x\rightarrow y$ or $x\nrightarrow y$, respectively. The notation $x\leftrightarrow y$ means that $x\rightarrow y$ and  $y\rightarrow x$.   For disjoint subsets $X$ and  $B$ of $V(D)$,  the notation  $X\rightarrow B$ means that every vertex  of $X$  dominates every vertex of $B$. $X\mapsto B$ means that $X\rightarrow B$ and there is no arc from a vertex of $B$ to a vertex of $X$.
   If $X\rightarrow B$ and $B\rightarrow F$, where $F$ is a subset of $V(D)$, we write $X\rightarrow B \rightarrow F$ (for short).
If $u\in V(D)$
   and $X=\{u\}$, we write $u$ instead of $\{u\}$. 
 The {\it out-neighborhood} of a vertex $x$ is the set $N^+(u)=\{y\in V(D)\, |\, uy\in A(D)\}$ and $N^-(u)=\{y\in V(D)\, |\, yu\in A(D)\}$ is the {\it in-neighborhood} of $u$. Similarly, if $X\subseteq V(D)$, then $N^+(u,X)=\{y\in X \,| \, uy\in A(D)\}$ and $N^-(u,X)=\{y\in A \, |\, yu\in A(D)\}$. 
The {\it out-degree} of $u$ is $d^+(u)=|N^+(u)|$ and $d^-(u)=|N^-(u)|$ is the {\it in-degree} of $u$. Similarly, $d^+(u,X)=|N^+(u,X)|$ and $d^-(u,X)=|N^-(u,X)|$. The {\it degree} of the vertex $u$ in $D$ is defined as $d(u)=d^+(u)+d^-(u)$ (similarly, $d(u,X)=d^+(u,X)+d^-(u,X)$). 

The subdigraph of $D$ induced by a subset $X$ of $V(D)$ is denoted by $D\langle X\rangle$.  
The  path (respectively, the cycle) consisting of the distinct vertices $x_1,x_2,\ldots ,x_m$ ($m\geq 2 $) of $V(D)$ and the arcs $x_ix_{i+1}$, $1\leq i\leq m-1$  (respectively, $x_ix_{i+1}$, $1\leq i\leq m-1$, and $x_mx_1$) of $A(D)$, is denoted by $x_1x_2\cdots x_m$ (respectively, $x_1x_2\cdots x_mx_1$). We say that
$x_1x_2\cdots x_m$  is an $(x_1,x_m)$-{\it path} in $D$, and if $m\geq 2$ and $x_1x_m\in A(D)$, then we say that $x_1x_2\cdots x_m$  is an $(x_1,x_m)$-{\it bypass of order} $m$ in $D$, in particular, if $m=|V(D)|$, then we have a {\it Hamiltonian bypass}. The {\it length} of a cycle or a path is the number of its arcs. A cycle of length $k$, $k\geq 2$, is denoted by $C_k$. A cycle passing through a vertex $x$ is denoted by $C(x)$. For a cycle  $C_k:=x_1x_2\cdots x_kx_1$, the subscripts considered modulo $k$, i.e., $x_i=x_s$ for every $s$ and $i$ such that  $i\equiv s\, (\hbox {mod} \,k)$. If $1\leq i\leq k$, then the {\it predecessor} of $x_i$ on $C_k$ is the vertex $x_{i-1}$ and is also denoted by $x^-_i$, the {\it successor} of $x_i$ on $C_k$ is the vertex $x_{i+1}$ and also denoted by $x^+_i$.  If $P$ is a path containing a subpath from $x$ to $y$, we let $P[x,y]$ denote that subpath. Similarly, if $C$ is a cycle containing vertices $x$ and $y$, $C[x,y]$ denotes the subpath of $C$ from $x$ to $y$.

A digraph $D$ is {\it strongly connected} (or just {\it strong}) if for every pair of distinct vertices $x,y$ there exists a path from $x$ to $y$ and a path from $y$ to $x$.  A digraph $D$ is $k$-{\it strongly connected} (or, $k$-{\it strong}), if $|V(D)|\geq k+1$ and  $D\langle V(D)\setminus X\rangle$ is strong for any set $X$ of at most $k-1$ vertices. 
   Two distinct vertices $x$ and $y$ of a digraph $D$ are {\it adjacent} if $xy\in A(D)$ or $yx\in A(D)$ (or both). Let $P=x_1x_2\ldots x_k$ and $Q=y_1y_2\ldots y_m$ be two vertex-disjoint paths in a digraph $D$, where $k\geq 1$ and $m\geq 1$. If $x_k\rightarrow y_1$, then  $PQ$  denotes the path $x_1x_2\ldots x_ky_1y_2\ldots y_m$. Let $m$ and $n$, $m\leq n$, be integers. By $[m,n]$ we denote the set $\{m,m+1,\ldots , n\}$.  The converse digraph of a digraph $D$ is the digraph  obtained from $D$ by reversing all arcs of $D$. We will use the {\it principle of digraph duality}: Let $D$ be a digraph, then $D$ contains a subdigraph $H$ if and only if the converse digraph of $D$ contains the converse subdigraph of $H$.

\section {Preliminaries}

The following well-known simple Lemmas 3.1-3.3 are the basis of our results
and other theorems about directed cycles and paths in digraphs. They are
will be used extensively in the proof of our result. 

\textbf{Lemma 3.1.} (H\"{a}ggkvist and Thomassen \cite{[16]}). {\it Let $D$ be a digraph of order  $p\geq 3$
 containing a
 cycle $C_m$, $2\leq m\leq p-1$. Let $x$ be a vertex not contained in this cycle. If  $d(x,V(C_m))\geq m+1$,
 then  for every $k$, $2\leq k\leq m+1$, $D$ contains a cycle of length $k$ that includes $x$.}

Lemma 3.2 is a  modification of a lemma by Bondy and Thomassen \cite{[7]}. 

\textbf{Lemma 3.2.} {\it Let $D$ be a digraph of order $p\geq 3$ containing a
 path $P:=x_1x_2\ldots x_m$, $2\leq m\leq p-1$ and  $x$ be a vertex not contained in that path.
  If one of the following conditions is true:

 $(i)$ $d(x,V(P))\geq m+2$;

 $(ii)$ $d(x,V(P))\geq m+1$ and $xx_1\notin A(D)$ or $x_mx\notin A(D)$;

 $(iii)$ $d(x,V(P))\geq m$, $xx_1\notin A(D)$ and $x_mx\notin A(D)$};

\noindent\textbf{} {\it then there is an $i$, $1\leq i\leq m-1$, such that
$x_ix,xx_{i+1}\in A(D)$,  i.e.,  $x_1x_2\ldots x_ixx_{i+1}\ldots x_m$ is a path of length $m$ in $D$ (we say that  $x$ can be inserted into $P$).}\\

The following lemma is a simple extension of a lemma by Bang-Jensen, Gutin and Li \cite{[3]}. 

\textbf{Lemma 3.3.} {\it Let $P=u_1u_2\ldots u_s$ be a path in a digraph $D$ (possibly, $s=1$) and let
$Q=v_1v_2\ldots v_t$ be a path (or $Q=v_1v_2\ldots v_tv_1$ be a cycle) in $D\langle V(D)\setminus V(P)\rangle$, $t\geq 2$. Suppose that for each $u_i$, $1\leq i\leq s$, there is an arc $v_jv_{j+1}$ on $Q$ such that $v_ju_i, u_iv_{j+1}\in A(D)$.
Then there exists a $(v_1,v_t)$-path (or a cycle) of length $t+k-1$ (respectively, $t+k$), $1\leq k\leq s$, with vertex set $\{v_1,v_2, \ldots , v_t\}\cup \{u_1,u_2,\ldots , u_k\}$}.\\

To prove Theorem 1.5, we also need the following lemma.

\textbf{Lemma 3.4} (Darbinyan \cite{[11]}). {\it Let $D$ be a strong digraph of order $p\geq 3$. If $p-1$ vertices of $V(D)$ have degree at least $p$, then $D$ is Hamiltonian or contains a cycle of length $p-1$  that contains all the vertices with degree at least $p$.}\\

Note that Lemma 3.4 is a direct consequence of Theorem 3.5.

\textbf{Theorem 3.5} (Berman and Liu \cite{[6]}). {\it Let $D$ be a strong digraph of order $p\geq 3$. Suppose that $M$ be a subset of $V(D)$ such that $d(x)+d(y)\geq 2p-1$ for every pair $x$, $y$ of distinct vertices in $M$ which are nonadjacent in $D$.  Then $D$   contains a cycle   that contains all the vertices of $M$.}\\

It is also worth mentioning that in \cite{[17]}, 
Li, Flandrin and Shu improved  Theorem 3.5 as follows: 

\textbf{Theorem 3.6} (Li, Flandrin and Shu \cite{[17]}). {\it Let $D$ be a  digraph of order $p\geq 3$ and let $S$ be a subset of $V(D)$ such that
 for every pair $x$, $y$ of distinct vertices of $S$ there is an $(x,y)$-path and a $(y,x)$-path in $D$. If   $d(x)+d(y)\geq 2p-1$ for any two nonadjacent distinct vertices $x\in S$ and $y\in S$,  then $D$   contains a cycle  through all the vertices of $S$.}\\

\section {Proofs of the main results}

To prove Theorem 1.5, we need the following Lemmas 4.1-4.4. The proofs of Lemmas 4.1, 4.2 and 4.3 can be found in \cite{[9]}.

   \textbf{Lemma 4.1} (Claim 1 and Case 1 in \cite {[9]}). {\it Let $D$ be a  digraph of order $p\geq 5$ and $z$ be some vertex of $D$. Suppose that every vertex  of  $V(D)\setminus \{z\}$ has degree   at least  $p$.  
 If $D$ contains a cycle of length at least $p-2$ through $z$, then $D$ contains a Hamiltonian bypass}.\\ 

Note that for $p=4$ Lemma 4.1 is not true. To see this, we consider a digraph of order four obtained from the complete digraph of order three and the complete digraph of order two by identifying a vertex of the first with a vertex of the second.\\
 
\textbf{Lemma 4.2} (Claim 4 and Proposition 1 in \cite{[9]}). {\it Let $D$ be a 2-strong digraph of order $p\geq 3$ and  $z$ be some vertex of $D$.   Suppose that every vertex  of  $V(D)\setminus \{z\}$ has degree   at least  $p$.
  If a longest path through $z$, say $P=x_1x_2\ldots x_{l}$, has length at least $p-3$ and $x_1x_l\in A(D)$, then  $D$ contains a Hamiltonian bypass}.\\
 
 Now, using Lemma 4.2, we will provide a new, significantly shorter proof of Lemma 4.3.
than that given in \cite{[9]}.

 \textbf{Lemma 4.3} (Lemma 4 in \cite {[9]}). {\it Let $D$ be a 2-strong digraph of order $p\geq 3$ and let $z$ be some vertex of $V(D)$. Suppose that every vertex of $V(D)\setminus \{z\}$ has degree  at least $p$. If a longest  cycle through $z$ has length $p-3$,   then $D$ contains a Hamiltonian bypass.}
 
 \textbf{ Proof}. Suppose, on the contrary, that there is a 2-strong digraph of order $p$ which satisfies the suppositions of Lemma 4.3, but has no Hamiltonian bypass. Let $C:=x_1x_2\ldots x_{p-3}x_1$  be a $C(z)$-cycle of length $p-3$ in $D$. Put $B:=V(D)\setminus V(C)=\{y_1,y_2,y_3\}$. Note that $d(y_i)\geq p$ since $z\in V(C)$. Since $C$ is a longest $C(z)$-cycle, it follows that every vertx $y_i$ cannot be inserted into $C$. 
 Then by Lemma 3.1,  for all $i\in [1,3]$, we have $d(y_i,V(C))\leq p-3$ and hence, $d(y_i,B)\geq 3$.  It is easy to see that  $D\langle B\rangle$ is strong. By the Ghouila-Houri theorem  $D\langle B\rangle$ is  Hamiltonian. Let $H:=y_1y_2y_3y_1$ be a Hamiltonian cycle in $D\langle B\rangle$. Using Lemma 4.2, we see that for any pair $i$, $j$ of the integers $i\in [1,p-3]$ and $j\in [1,3]$, $d^-(y_j,\{x_i,x_{i+1}\})\leq 1$ and $d^+(y_j,\{x_i,x_{i+1}\})\leq 1$. We also have, if $x_i\rightarrow y_j$, 
 then $y_jx_{i+1}\notin A(D)$ since $C$ is a longest $C(z)$-cycle. 
  Using this and Lemma 4.2, it is not difficult to see that if 
 $x_i\rightarrow y_j$   (respectively, $y_j\rightarrow x_i$), then 
 $d^+(x_{i+1},B)=0$ (respectively, $d^-(x_{i+1}, B)=0$). 
 Since $D$ is 2-strong, there are two distinct vertices $x_l$ and $x_k$ such that  
 $d^+(x_{l},B)\geq 1$ and  $d^+(x_{k},B)\geq 1$. By the above arguments, it is easy to see that
  $x_l\notin \{x_{k-1},x_{k+1}\}$. Therefore, we can choose a vertex $x_k$ such that $d^+(x_{k},B)\geq 1$, say $x_k\rightarrow y_1$, and $z\notin \{x_{k+1},x_{k+2}\}$. Then it is not difficult to see that $d(x_{k+1},B)=0$, i.e.,
  $d(x_{k+1})=d(x_{k+1},V(C[x_{k+2},x_k]))\geq p$.
  Then by Lemma 3.2, there exsits  an $(x_{k+2},x_k)$-path, say $R$, with vertex set  $V(C)$. This together with Lemmas 4.1 and 4.2 implies that $d(x_{k+2},B)=0$, in particular, we have that $d(y_3,\{x_{k+1},x_{k+2}\})=0$. Therefore,
   $d(y_3,V(C[x_{k+3},x_k]))\geq p-4$. Since $y_3$ cannot be inserted into $C$ and  $|V(C[x_{k+3},x_k])|=p-5$, using Lemma 3.2, we obtain $y_3\rightarrow x_{k+3}$. This means that $C[x_{k+3},x_k]y_1y_2y_3x_{k+3}$ is a $C(z)$-cycle of length $p-2$, a contradiction. Lemma 4.3 is proved. \fbox \\\\  
 
 \textbf{Lemma 4.4.} {\it Let $D$ be a 2-strong digraph of order $p\geq 3$ and let $z$ be some vertex of $D$.  Suppose that every vertex of $V(D)\setminus \{z\}$ has degree  at least $p$. If a longest cycle through $z$ has length $p-4$,   then $D$ contains a Hamiltonian bypass.}
 
 \textbf{Proof}. Suppose, on the contrary, that there exists a 2-strong digraph of order $p$ which satisfies the assumptions of Lemma 4.4, but has no Hamiltonian bypass. Let $C:=x_1x_2\ldots x_{p-4}x_1$  be a cycle of length $p-4$ through $z$ in $D$. Put $B:=V(D)\setminus V(C)=\{y_1,y_2,y_3,y_4\}$. Since $C$ is a longest $C(z)$-cycle, it follows that every vertx $y_i$ cannot be inserted into $C$. Then by Lemma 3.1,  for all $i\in [1,4]$ we have $d(y_i,V(C))\leq p-4$ and hence, $d(y_i,B)\geq 4$.

 Now, we divide the proof into two cases depending on  whether
  $D\langle B\rangle$ is strong or not.
 
 \textbf{Case 1.}  $D\langle B\rangle$ is not strong.
 
 Since $|B|=4$, $D\langle B\rangle$ is not strong and for all $y_i\in B$, $d(y_i,B)\geq 4$, it follows that $d(y_i,B)= 4$ and $d(y_i,V(C))= p-4$. Without loss of generality, we can assume that $y_1\leftrightarrow y_2$, $y_3\leftrightarrow y_4$ and $\{y_1,y_2\}\mapsto \{y_3,y_4\}$.
 
 \textbf{Subcase 1.1.} There are vertices $y\in B$ and $u\in V(C)$ which are not adjacent.
 
 By the digraph duality and symmetry between $y_1$ and $y_2$ ($y_3$ and $y_4$), we can assume that $y=y_1$. We know that $d(y_1,V(C[u^+,u^-]))=p-4$.
 Since $y_1$ cannot be inserted into $C[u^+,u^-]$ and $|V(C[u^+,u^-])|=p-5$, using Lemma 3.2, we obtain $u^-\rightarrow y_1 \rightarrow u^+$. 
  Since $C$ is a longest $C(z)$-cycle, it is easy to see that $d^-(u, B)=0$. If $d^+(u, B)\not= 0$, 
   then  $d^+(u,\{y_1,y_2\})= 0$ and hence, for some $i\in [3,4]$, $u\rightarrow y_i$. Therefore,   $C[u,u^-]y_1y_2y_i$ is a $(u,y_i)$-bypass of order  $p-1$, and by Lemma 4.2, $D$ has a Hamiltonian bypass, a contradiction. 
   We may therefore assume that $d(u,B)=0$. This and Lemma 3.2 imply that $u^-\rightarrow B\rightarrow u^+$ since $C$ is a longest $C(z)$-cycle and for all $y\in B$, $d(y,V(C))=p-4$. Now we consider the cycle $Q:=y_1C[u^+,u^-]y_1$ of length $p-4$. We see that $u$ cannot be inserted into $Q$ since $C$ is a longest $C(z)$-cycle. Then by Lemma 3.1, $d(u,V(Q))\leq p-4$. This together with $d(u,B)=0$ implies that $u=z$. Using this and Lemma 3.2,
   it is not difficult to show that every pair of vertices $y\in B$ and $v\in V(C)\setminus \{z\}$ are adjacent. 
   For the  sake of clarity, let $u=x_{p-4}$. Then $u^-=x_{p-5}$ and $u^+=x_1$. 
   Since $C$ is a longest $C(z)$-cycle and $y_1$ ($y_2$) is adjacent to every vertex of $V(C)\setminus \{x_{p-4}=z\}$, it follows that there is an integer $l\in [1,p-5]$ such that $\{x_l,x_{l+1},\ldots , x_{p-5}\}\rightarrow \{y_1,y_2\}\rightarrow \{x_1,x_2,\ldots , x_l\}$. It is easy to see that $l\leq p-6$ (for otherwise, if $l=p-5$, then in $D-x_{p-5}$ there is no arc from a vertex of $V(D)\setminus \{x_{p-5},y_1,y_2\}$ to a vertex of $\{y_1,y_2\}$, which contradicts that $D$ is 2-strong). Then $C[x_{p-5},x_{p-6}]y_2y_1$ is an $(x_{p-5},y_1)$-bypass of order $p-2$. Therefore by Lemma 4.2, $D$ contains a Hamiltonian bypass, a contradiction.

 \textbf{Subcase 1.2.} Every pair $y$, $u$ of vertices $y\in B$ and $u\in V(C)$ are adjacent.
 
 Without loss of generality, we consider the vertex $y_1$. By the maximality of the cycle $C$, either $V(C)\mapsto y_1$ or $y_1\mapsto V(C)$. If $V(C)\mapsto y_1$, then obviously, $A(B\rightarrow V(C))=\emptyset$, which contradicts that $D$ is strong. So we can   assume that $y_1\mapsto V(C)$. Then $y_2\mapsto V(C)$. Hence, $y_2y_1C[x_2,x_1]$ is a $(y_2,x_1)$-bypass through $z$ of order $p-2$. Then by Lemma 4.2, $D$ has a Hamiltonian bypass, a contradiction.  This contradiction completes the discussion of Case 1.
 
 \textbf{Case 2.}  $D\langle B\rangle$ is 
  strong.
  
  Since $d(y_i,B)\geq 4$ for all $y_i\in B$ by the Ghouila-Houri theorem, $D\langle B\rangle$ is Hamiltonian. Let $H:=y_1y_2y_3y_4y_1$ be a Hamiltonian cycle in $D\langle B\rangle$. For this case, we first prove the following claim.
  
 \textbf{Claim 1}.
 
  (i). {\it Suppose that $x_k\rightarrow y_i$, where $k\in [1,p-4]$ and $i\in [1,4]$. Then the following holds:
  
  (a) $d^-(x_{k+1},B)=0$ and (b) $d^+(x_{k+1},B\setminus \{y_i\})=0$.  In particular,  $d(x_{k+1},B\setminus \{y_i\})=0$.
  
  (ii) Suppose that $y_i\rightarrow x_k$, where $k\in [1,p-4]$ and $i\in [1,4]$. Then the following holds: 
  
  (a) $d^+(x_{k-1},B)=0$ and  (b) $d^-(x_{k-1},B\setminus \{y_i\})=0$. In particular,  $d(x_{k-1},B\setminus \{y_i\})=0$.
  
  (iii) For each pair of integers  $k\in [1,p-4]$ and $i\in [1,4]$,  the following hold: $d^-(y_i,\{x_k,x_{k+1}\})\leq 1$ and $d^+(y_i,\{x_k,x_{k+1}\})\leq 1$}. In particular, $d(y_i,\{x_k,x_{k+1}\})\leq 2$.
  
  \textbf{Proof.} For the sake of clarity, let  $y_i=y_1$.
  
  i(a). Suppose, on the contrary, that is for some $k\in [1,p-4]$,   $x_k\rightarrow y_1$ and   $d^-(x_{k+1},B)\not=0$, i.e.,         $y_j\rightarrow x_{k+1}$ with $j\in [1,4]$. Then $C[x_{k+1},x_k]H[y_1,y_j]x_{k+1}$ is a $C(z)$-cycle of length at least $p-3$, which contradicts that a longest $C(z)$ cycle has length $p-4$.
  
  i(b).  Suppose, on the contrary, that is for some $k\in [1,p-4]$,    $x_k\rightarrow y_1$ and  $d^+(x_{k+1},B\setminus \{y_1\})\not=0$. Then
  $x_{k+1}\rightarrow y_j$ with $j\in[2,4]$ and $C[x_{k+1},x_{k}]H[y_1,y_j]$ is an
   $(x_{k+1},y_j)$-bypass through $z$ of order at least $p-2$. Then by Lemma 4.2, $D$ has a Hamiltonian bypass, a contradiction.  

   From i(a) and i(b) it follows that if $x_k\rightarrow y_j$ with $k\in [1,p-4]$ and $j\in [1,4]$, then $d(x_{k+1},B\setminus \{y_j\})=0$.

  (ii). From the digraph duality and Claim 1(i)  it follows that Claim 1(ii) is also true.
  
  (iii). Assume for a contradiction, that for some $y_i$, say $y_1$, and $x_k$ with $k\in [1,p-4]$, 
  $\{x_k,x_{k+1}\}\rightarrow y_1$ or $y_1\rightarrow \{x_k,x_{k+1}\}$. By the digraph duality, it  suffices  to consider only the case   $\{x_k,x_{k+1}\}\rightarrow y_1$. Then by Claim 1(i), we have
  $$
  d^+(x_k,B\setminus \{y_1\})=d^+(x_{k-1},B\setminus \{y_1\})=0. 
                      \eqno (1)
  $$ 
 Again using Claim 1(i), for each $j\in [2,4]$ we obtain,  $d(y_j,\{x_{k+1},x_{k+2}\})=0$.  This together with $d(y_j)\geq p$ implies that  $d(y_j,V(C[x_{k+3},x_k]))\geq p-6$. Then, since $x_ky_j\notin A(D)$ and  $y_j$ cannot be inserted into 
$C[x_{k+3},x_k]$, by Lemma 3.2 we have that
   $d(y_j,V(C[x_{k+3},x_k]))
  = p-6$, $d(y_j,B)=6$ and $\{y_2,y_3,y_4\}\rightarrow x_{k+3}$, in particular, $D\langle B\rangle$ is a complete digraph. It is easy to see that $x_{k+2}=z$ and $d(z,B)=0$. From (1) and $d(y_j,V(C[x_{k+3},x_k]))= p-6$ it follows that
  $\{y_2,y_3,y_4\}\rightarrow x_k$. 
  Assume that $x_{k-2}y_j\in A(D)$ with $j\in [1,4]$. Then, since $x_{k+2}=z$, for $j\not=2$ we have $C[x_k,x_{k-2}]H[y_j,y_{j+3}]x_k$ is a $C(z)$-cycle of length $p-1$, and if $j=2$, then $C[x_k,x_{k-2}]H[y_j,y_{j+2}]x_k$ is a $C(z)$-cycle of length $p-2$.  Thus in bout cases we get a contradiction. Therefore, 
    $d^+(x_{k-2},B)=0$. Then by Lemma 3.2, for each $j\in [2,4]$ we have
\begin{eqnarray*}
  p &\leq & d(y_j) =d(y_j,B)+d(y_j,\{x_k\}) \\
  & + & d(y_j,V(C[x_{k+3},x_{k-2}])) \\
  &=& 6+1+p-8=p-1,
 \end{eqnarray*} 
   a contradiction.
   This proves that for all  $i\in [1,4]$ and $k\in [1,p-4]$,  $d^-(y_i,\{x_k,x_{k+1}\})\leq 1$.
      Claim 1 is proved. \fbox \\\\

  Since $D$ is 2-strong, there are two distinct vertices $x_s$, $x_t\in V(C)$ such that $d^+(x_s,B)\geq 1$ and  $d^+(x_t,B)\geq 1$. From Claim 1 it follows  that
  $x_s\notin \{x_{t-1},x_{t+1}\}$. Therefore, we can choose a vertex $x_k$ such that $d^+(x_k,B)\geq 1$ and  $z\notin \{x_{k+1},x_{k+2}\}$. Without loss of generality assume that  $x_k\rightarrow y_1$. By Claims 1(i) and 1(iii), 
  $$ 
  d(x_{k+1},B)= d^+(x_{k-1},B)=0.  \eqno (2) 
  $$ 
  
Since $z\notin \{x_{k+1},x_{k+2}\}$, we have that  $d(x_{k+1})\geq p$ and $d(x_{k+2})\geq p$.  From $d(x_{k+1},B)=0$
it follows that
 $$
 p\leq d(x_{k+1})= d(x_{k+1},V(C[x_{k+2},x_k]))\leq 2(p-5),
 $$ 
 in particular, $p\geq 10$.
    By Lemma 3.2, there exists an $(x_{k+2},x_k)$-path, say $Q$, with vertex set $V(C)$. Therefore, if  $y_i\rightarrow x_{k+2}$ with $i\in [1,4]$, then  $H[y_1,y_i]Qy_1$ is a $C(z)$-cycle of length at least $p-3$, a contradiction. We may therefore assume that $d^-(x_{k+2},B)=0$. If $x_{k+2}\rightarrow y_i$ with $i\not=1$, then $QH[y_1,y_i]$ is an $(x_{k+2},y_i)$-bypass of order at least $p-2$. Therefore by Lemma 4.2, $D$ contains a Hamiltonian bypass, a contradiction.  Thus, we may assume that 
  $$
  d(x_{k+2},\{y_2,y_3,y_4\})=0. \eqno (3)
  $$ 
  This and $d(x_{k+1},B)=0$ (by (2)) imply that there is no arc between $\{x_{k+1},x_{k+2}\}$ and 
$\{y_{2},y_{3},y_4\}$.

Now we consider two subcases depending on whether $y_1x_{k-1}\in A(D)$ or not.

 \textbf{Subcase 2.1.} $y_1 x_{k-1}\in A(D)$.
 
 Since $C$ is a longest $C(z)$-cycle, we have $d^+(x_{k-2},B)=0$. On the other hand, Claims 1 implies that $d^-(x_{k-2},B)=0$. Then 
 $d(x_{k-2},B)=0$ and  for every $i\in [2,4]$,  $d(y_i,\{x_{k-2},x_{k+1},x_{k+2}\})=0$. Using the fact that $y_1\rightarrow x_{k-1}$ and Claim 1(ii), we obtain $d^-(x_k,\{y_2,y_3,y_4\})=0$. Since every vertex $y_i$ cannot be inserted into $C[x_{k+3},x_{k-3}]$ and $d(y_i)\geq p$, using Lemma 3.2 and the equality $d^+(x_{k-1},B)=0$ (by (2)), for every $i\in [2,4]$ we obtain,
\begin{eqnarray*}
  p &\leq & d(y_i) =d(y_i,B)+d(y_i,V(C[x_{k+3},x_{k-3}])) \\
  & + & d^+(y_i,\{x_{k-1}\})+d^-(y_i,\{x_{k}\})  \\
  &\leq& 6+p-8+2=p,
 \end{eqnarray*}

This implies that for each $i\in [2,4]$,
  $d(y_i,V(C[x_{k+3},x_{k-3}]))=p-8$ and   by Lemma 3.2, $x_{k-3}\rightarrow \{y_2,y_3,y_4\}\rightarrow x_{k+3}$. So $x_ky_1y_2y_3y_4C[x_{k+3},x_k]$ is a $C(z)$-cycle of length $p-2$, a contradiction.

 \textbf{Subcase 2.2.} $y_1x_{k-1}\notin A(D)$.
 
 Then by Claim 1(iii), the vertices $y_1$ and $x_{k-1}$ are not adjacent since $x_k\rightarrow y_1$. By $d(x_{k+2})\geq p$, $y_1x_{k+2}\notin A(D)$ and $d(x_{k+2},\{y_2,y_3,y_4\})=0$ we have
\begin{eqnarray*}
  p &\leq & d(x_{k+2}) = d(x_{k+2},\{y_2,y_3,y_4\})
   \\
  & + & d(x_{k+2},\{y_1,x_{k+1}\})+  d(x_{k+2},V(C[x_{k+3},x_{k}]))   \\
  &\leq& d(x_{k+2},V(C[x_{k+3},x_{k}]))+3.
 \end{eqnarray*}

 Therefore,
 $d(x_{k+2},V(C[x_{k+3},x_{k}]))\geq p-3$. 
 Therefore by Lemma 3.2, $x_{k+2}$  can be inserted into $C[x_{k+3},x_{k}]$ since $|V(C[x_{k+3},x_{k}])|=p-6$. So each of the vertices $x_{k+1}$ and $x_{k+2}$ can be inserted into $C[x_{k+3},x_{k}]$. Therefore,
  there exists an $(x_{k+3},x_k)$-path , say $R$, with vertex set $V(C)$. Since $C$ is a longest $C(z)$-cycle, it is not difficult to see that 
 $$
 d^-(x_{k+3},B)=d^+(x_{k+3},\{y_2,y_3,y_4\})=0 \eqno (4)   
 $$ 
 Indeed, if $y_i\rightarrow x_{k+3}$, then $y_iRH[y_1,y_i]$ is a $C(z)$-cycle of length at least $p-3$, if $x_{k+3}\rightarrow y_j$ with $j\not= 1$, then $RH[y_1,y_j]$ is an $(x_{k+3},y_j)$-bypass of order at least $p-2$ and by Lemma 4.2, $D$ has a Hamiltonian bypass. So, in both cases we have a contradiction, which shows that (3) is true. (4) together with (3) and (2) implies that for every $i\in [2,4]$, 
 $d(y_i,\{x_{k+1},x_{k+2},x_{k+3}\})=0$. Then $d(y_i,V(C[x_{k+4},x_k]))\geq p-6$ and  by Lemma 3.2, $x_k\rightarrow \{y_2,y_3,y_4\}\rightarrow x_{k+4}$ and 
 $x_ky_1y_2y_3y_4x_{k+4}x_{k+5}\ldots x_k$ is a cycle of length $p-3$. 
 Therefore, $z=x_{k+3}$. Note that $z\not=x_{k-1}$ since $p\geq 10$.  
 Then, $d(x_{k-1})\geq p$ 
 and 
 $d(x_{k-1},V(C[x_k,x_{k-2}])\geq p-3$ 
 since $d(x_{k-1},\{y_1\})=d^+(x_{k-1},\{y_2,y_3,y_4\})=0$.
 Using Lemma 3.2 and $|V(C[x_k,x_{k-2})|=p-5$, we conclude that there exists  an $(x_k,x_{k-2})$-path, say $F$, with vertex set  $V(C)$.
  If $x_{k-2}\rightarrow y_i$ with $i\not= 1$, then  the path $FH[y_i,y_1]$
   is an $(x_k,y_i)$-bypass of order at least $p-2$, a contradiction.  We may therefore assume that $d^+(x_{k-2}, \{y_2,y_3,y_4\})=0$.
    Using Claim 1(ii), it is not difficult to show that $d(y_i,\{x_{k-1},x_k\})\leq 2$. Then, since every vertex $y_i$ with $i\in [2,4]$ cannot be inserted into $C[x_{k+4},x_{k-1}]$,
   using Lemma 3.2, we obtain 
\begin{eqnarray*}
   p &\leq & d(y_i)=d(y_i,B)+d(y_i,\{x_{k-1})
  + d(y_i,V(C[x_{k+4},x_{k-2}])) 
   \\   
  &\leq &  6+2+p-9=p-1,
 \end{eqnarray*} 
 which is a contradiction. This contradiction completes the proof of Lemma 4.4. \fbox
 \\\\
 
 We are now ready to prove our main results. \\

{\it The proof of Theorem 1.5}. Suppose, on the contrary, that is $D$ satisfies the conditions of the theorem but does not contain a Hamiltonian bypass. Obviously, $p\geq 5$. From Lemma 4.4 it follows that $D$  does not contain a $C(z)$-cycle of length greater than $p-5$. Since $D$ is not Hamiltonian, by Lemma 1.1 it contains a cycle of length $p-1$. 
Let $R:=x_1x_2\ldots x_{p-1}x_1$ be a cycle of length $p-1$ in $D$. Then by Lemmas 1.1 and 4.4, we know that $z\notin V(R)$ and for every $i$, $j$, $1\leq i, j\leq p-1$ the following hold:
 $$
 d^-(z,\{x_i\})+d^+(z,\{x_{i+1},x_{i+2},x_{i+3},x_{i+4},x_{i+5}\})\leq 1,
 $$
 $$
 d^+(z,\{x_j\})+d^-(z,\{x_{j-1},x_{j-2},x_{j-3},x_{j-4},x_{j-5}\})\leq 1.
 $$
 Therefore, $d^-(z)+5d^+(z)\leq p-1$ and $5d^-(z)+d^+(z)\leq p-1$. This means that $6d(z)\leq 2p-2$, i.e., $d(z)\leq (p-1)/3$, which contradicts that $d(z)>(p-1)/3$. The theorem is proved. \fbox \\\\
 
 \textbf{Corollary 1} (Benhocine \cite{[4]}). {\it Every strong digraph $D$ of order $p\geq 3$ and with minimum degree at least $p$ contains $D(p,2)$}.
 
 \textbf{\it Proof:}  According to the Ghoula-Houri theorem, $D$ is Hamiltonian. Therefore, it follows from Lemma 4.1 that $D$ contains a Hamiltonian bypass. \fbox \\
 
 \textbf{Corollary 2} (Darbinyan \cite{[9]}). {\it Let $D$ be a 2-strong digraph of order $p\geq 3$ and let  $z$ be some vertex of $D$. Suppose that $d(x)\geq p$ for each vertex $x\in  V(D)\setminus \{z\}$.  If $D$ is Hamiltonian or $d(z)>0.4(p-1)$,  then $D$ has a Hamiltonian bypass.}\\

We prove the following theorem for a balanced bipartite digraph.

\textbf{Theorem 4.1.} {\it Let $B(X,Y)$ be
a balanced bipartite digraph with partite sets $X=\{x_1,x_2,\ldots ,\\ x_a\}$ and $Y=\{y_1,y_2,\ldots , y_a\}$ satisfying the condition $A_l$, where $l\geq 0$. Suppose that $B(X,Y)$ has a Hamiltonian cycle $C=x_1y_1x_2y_2\ldots x_ay_ax_1$ such that $y_{i+k-1}x_i\in A(B(X,Y))$ (or $x_{i+k-1}y_i\in A(B(X,Y))$), where $k$ is the minimum possible with this properties. If $B(X,Y)$ has no Hamiltonian bypass, then  $l=0$ and $k=2$}. 

 \textbf{Proof.} Let $B(X,Y)$ be a balanced bipartite digraph  satisfying the conditions of the theorem.  Without loss of generality, we may assume that $i=1$, i.e., $y_kx_1\in A(B(X,Y))$. Then it is easy to see that $k\geq 2$ and 
 $$
 d^+(x_k,\{y_1,y_2,\ldots , y_{k-1}\})=d^-(y_1,\{x_2,x_3,\ldots , x_{k}\})=0.
 $$
If for some $j\in [k+1,a]$, $x_jy_1\in A(B(X,Y))$ and $x_ky_j\in A(B(X,Y))$, then $y_kx_{k+1}y_{k+1}\ldots \\x_jy_1x_2y_3\ldots x_ky_jx_{j+1}\ldots y_ax_1$ is a Hamiltonian bypass in $B(X,Y)$ since  $y_kx_1\in A(B(X,Y))$, which is a contradiction. We may therefore assume that for all $j\in [k+1,a]$, $$d^-(y_1,\{x_j\})+d^+(x_k, \{y_j\})\leq 1.$$ Then
\begin{eqnarray*}
  a+l &\leq & d^-(y_1)+d^+(x_k)=d^-(y_1,\{x_{k+1},x_{k+2},\ldots , x_{a}\})
   \\
  & + & d^+(x_k,\{y_{k+1},y_{k+2},\ldots , y_{a}\})+ d^-(y_1,\{x_1\})+d^+(x_k,\{y_k\}) \\ 
  & \leq & a-k+2.
 \end{eqnarray*} 
 
Therefore, $l\leq 2-k$ and $l=0$, $k=2$ since $l\geq 0$ and $k\geq 2$.  The theorem is proved. \fbox\\\\

 We will now prove Theorem 1.6 using Wang's theorem [21] and Theorem 4.1.

{\it The proof of theorem 1.6.} 
According to Wang's theorem, if B is different from B6, then
 B is Hamiltonian.
 Therefore by Theorem 4.1 and the fact that $B_6$ has a Hamiltonian bypass (see the end of Section 1), we have that $B$ has a Hamiltonian bypass. \fbox\\\\
  
The following example (see Figure 1(a), in Figgure 1(a) an undirected edge represents an arc with a clockwise direction) shows that there is a balanced bipartite Hamiltonian  digraph of order 8 satisfying the condition $A_0$, which has no Hamiltonian bypass. We don't know if there is a balanced bipartite strong digraph of order $2a\geq 10$  which satisfies the condition $A_0$ and has no Hamiltonian bypass. We propose the following problem.  

\textbf{Problem 1.}  Let $B$ be a balanced bipartite Hamiltonian  digraph of order $2a\geq 10$ satisfying the condition $A_0$. Decide when $B$ has a Hamiltonian bypass?

\section {Conclusion}

 In the current article, we investigated the existence of a Hamiltonian bypass in  2-strong digraphs of order $p$,  in which $p-1$ vertices have degrees at least $p$. We proved that if such digraphs are Hamiltonian or have the
minimum degree is greater than $(p-1)/3$, then they contain a Hamiltonian bypass. 
Now, let us   consider the digraph $H(n)$. As mentioned in Section 1,  $H(n)$ is 2-strong, $d(x_0)=4$ and is not Hamiltonian. It is not difficult check that  $y_1y_2y_3x_2x_3\ldots x_{n-4}x_0x_1$ is a Hamiltonian bypass in $H(n)$.
By this argument,
we believe that  the following conjecture is true. 

\textbf{Conjecture 1.}  {\it Let $D$ be a 2-strong digraph of order $p$. If $p-1$ vertices in $V(D)$ have degrees at least $p$, then $D$ contains a  Hamiltonian bypass.}\\

In \cite{[13]}, the following theorem was proved.

\textbf{Theorem 5.1.} {\it Let $D$ be a 2-strong digraph of order $p\geq 3$ such that for every distinct pair of nonadjacent vertices $x$, $y$ and $w$, $z$ we have $d(x)+d(y)+d(w)+d(z)\geq 4p-3$. Then $D$ is Hamiltonian.} \\

For the existence of Hamiltonian   bypasses in digraphs, we propose the following conjecture.

\textbf{Conjecture 2.} {\it Let $D$ be a strong digraph of order $p$ such that for every distinct pair of nonadjacent vertices $x$, $y$ and $w$, $z$ we have $d(x)+d(y)+d(w)+d(z)\geq 4p-4$. Then $D$ has a  Hamiltonian bypass, except for certain digraphs.} \\

Note that if in Conjecture 2, we assume that $d(u)+d(v)\geq 2p-2$ for all pairs of nonadjacent vertices $u$ and , $v$, then it follows from Theorem 1.3  that $D$ has a Hamiltonian bypass, except for certain digraphs. Therefore, we can assume that $D$ contains exactly one pair of nonadjacent distinct vertices $u$, $v$ such that  $d(u)+d(v)\leq 2p-2-k$, where $k\geq 1$ is an integer. In this case, we have $d(x)+d(y)\geq 2p-2+k$ for all pairs of nonadjacent distinct vertices $x$ and $y$, exsept for the pair   $u$ and  $v$.\\

For arbitrary integers $n\geq 3$ and $q\in [2,n]$, $D(n,q)$ denotes the digraph of order $n$ obtained from a directed cycle of length $n$, by changing the orientation of $q-1$ consecutive arcs.  Benhocine and  Wojda \cite{[5]} proved that any digraph $D$ of order $p$ that satisfies the condition that the sum of the degrees for any two non-adjacent vertices is at least $2p$, contains $D(n,2)$ for every $n\in [3,p]$, except for certain digraphs. 
We are sure that if we replace $2p$ with $2p-1$ in this result, then the result will be correct again (Conjecture 3). 

\textbf{Conjecture 3.} {\it Let $D$ be a strong digraph of order $p$ which satisfies the condition that the sum of the degrees for any two nonadjacent vertices is at least $2p-1$. Then $D$ contains $D(n,2)$ for each $n\in [3,p]$, except for certain digraphs}.

Note that, the author \cite{[8]} proved that,  under the conditions of Conjecture 3 a digraph contains a $D(p,3)$.\\

 Declaration of Competing Interest
 
The authors declare that they have no known competing financial interests or personal relationships that could have appeared to influence the work reported in this paper.

Data availability

No data was used for the research described in the article.

\textbf{Aknowledgements:}

 The author would like to thank Dr Parandzem Hakobyan for formatting the manuscript of this paper.


\begin{thebibliography}{25}




\bibitem{[1]} J. Adamus and  L. Adamus, A degree condition for cycles of maximum length in bipartite digraphs, {\it Discrete Math.}, 312(6) (2012), pp. 1117-1122, 10.1016/j.disc.2011.11.032

\bibitem{[2]} J. Bang-Jensen and G. Gutin, "Digraphs: Theory,  Algorithms and Applications", Springer, 2000.

\bibitem{[3]} J. Bang-Jensen, G. Gutin and H. Li, Sufficient conditions for a digraph to be Hamiltonian, {\it J. Graph Theory},  22(2) (1996),pp. 181-187, 
 10.1002/(SICI)1097-0118(199606)22:2 <181::AID-JGT9> 3.0.CO;2-J
 
 
 \bibitem{[4]} A. Benhocine, On the existence of a specified cycle in digraphs with constraints on degrees, {\it J. Graph Theory},  8 (1984), pp. 101-107, 10.1002/jgt.3190080111
 
 \bibitem{[5]} A. Benhocine and A.P. Wojda, Bypasses in Digraphs, 
 {\it Ars Combinatoria},  16 (1983), pp.  85-94

\bibitem{[6]} K.A. Berman and X. Liu, Cycles through Large Degree Vertices in Digraphs: A Generalization of Meyniel's Theorem, {\it J. Combin.Theory} B 74 (1998), pp. 20-27, 
10.1006/jctb.1998.1829
 
 \bibitem{[7]} J.A. Bondy and C. Thomassen, A short proof of Meyniel's theorem, {\it Discrete Math.},  19 (1977), pp. 195-197, 10.1016/0012-365X(77)90034-6
 
\bibitem{[8]} S.Kh. Darbinyan, On Hamiltonian bypasses in digraphs satisfying Meyniel-like condition, {\it Math. Problems of Computer Science},  
 20 (1998), pp. 7-19

\bibitem{[9]} S.Kh. Darbinyan, A Note on Hamiltonian Bypasses in digraphs with large degrees, {\it Math. Problems of Computer Science},  
 54 (2020), pp. 7-17,  10.51408/1963-0055


\bibitem{[10]} S.Kh. Darbinyan, On Hamiltonian bypasses in digraphs with the condition of Y. Manoussakis, {\it "2015  Computer Science and Information Technologies (CSIT)}, Yerevan,  (2015), pp. 53-63, 10.1109/CSITechnol.2015.7358250 

\bibitem{[11]} S.Kh. Darbinyan, On Three Conjectures of Thomassen and the Extremal
Digraphs for Two Conjectures of Nash‐Williams, {\it J. of Graph Theory}, 109,4 (2025), pp. 412-425, https://doi.org/10.1002/jgt.23233


\bibitem{[12]} S.Kh. Darbinyan, A new sufficient condition for a 2-strong digraph to be Hamiltonian, {\it Discrete Mathematics and Theoretical  Computer  Science}, vol. 26:2 (2024), no. 7, 10.46298/dmtcs.11560

\bibitem{[13]} S.Kh. Darbinyan, A new sufficient condition for a digraph to be Hamiltonian- A proof of Manoussakis conjecture, {\it Discrete Mathematics and Theoretical  Computer  Science}, 22(4) (2021), no. 12, 10.23638/DMTCS-22-4-12

\bibitem{[14]} S.Kh. Darbinyan and I.A. Karapetyan, On Hamiltonian bypasses in one class of Hamiltonian digraphs, {\it Math. Problems of Computer Science}, 41 (2014), pp. 23-37 

\bibitem{[15]} A.~Ghouila-Houri, 
Une condition suffisante d'existence d'un circuit hamiltonien.
 \emph{C. R. Acad. Sci. Paris Ser. A-B}, 251 (1960), pp. 495-497

\bibitem{[16]} R. H\"{a}ggkvist and C. Thomassen, On pancyclic digraphs, {\it Journal Combinatorial Theory ser. B},   20 (1976), pp. 20-40, 10.1016/0095-8956(76)90063-0

\bibitem{[17]} H. Li, E. Flandrin and J. Shu, A sufficient condition for cyclability in directed graphs, {\it Discrete Math.} 307 (2007), pp. 1291-1297, 10.1016/j.disc.2005.11.066

\bibitem{[18]}  Y. Manoussakis, Directed hamiltonian graphs, {\it J. Graph Theory}, 16(1) (1992), pp. 51-59, 10.1002/jgt.3190160106


\bibitem{[19]} M.~Meyniel,
  Une condition suffisante d'existence d'un circuit hamiltonien dans un   graphe oriente,
 \emph{J. Combin. Theory Ser. B}, 14  (1973), pp. 137-147,    
10.1016/0095-8956(73)90057-9

\bibitem{[20]} C.~St. J.~A. Nash-Williams, 
 Hamilton circuits in graphs and digraphs, the many facets of graph
  theory.
 \emph{Springer Verlag Lecture Notes}, 110 (1969), pp. 237-243

\bibitem{[21]} R. Wang, Extremal digraphs on Woodall-type condition for hamiltonian cycles in balanced bipartite digraphs,
{\it J. Graph Theory}, 97(2) (2021), pp. 185-207,  10.1002/jgt.22649

\bibitem{[22]} D. Woodall,
 Sufficient conditions for circuits in graphs.
 \emph{Proc. Lond. Math. Soc.}, 24 (1972), pp. 739--755,  10.1112/plms/s3-24.4.739 




\end{thebibliography}
\end{document}